\documentclass[10pt,psamsfonts]{amsart}
\usepackage{amsmath}
\usepackage{amsthm}
\usepackage{amssymb}
\usepackage{amscd}
\usepackage{amsfonts}
\usepackage{amsbsy}
\usepackage{graphicx}
\usepackage[dvips]{psfrag}
\usepackage{array}
\usepackage{color}
\usepackage{epsfig}
\usepackage{url}
\usepackage{overpic}
\usepackage{tikz-cd}
\usepackage{enumitem}
\usepackage{hyperref}
\usepackage{soul}
\usepackage{faktor}
\usepackage{float}
\usepackage[normalem]{ulem}
\usepackage{nccmath}
\usepackage[foot]{amsaddr}
\usepackage{cleveref}

\usepackage{textcomp}

\newcolumntype{L}{>{\displaystyle}l}
\newcolumntype{C}{>{\displaystyle}c}
\newcolumntype{R}{>{\displaystyle}r}

\newcommand{\R}{\ensuremath{\mathbb{R}}}

\def\p{\partial}
\def\e{\varepsilon}

\newtheorem {theorem} {Theorem}

\newtheorem {lemma}{Lemma}

\newtheorem{hypothesis}{Hypothesis}
\newtheorem{wassumption}{Working Assumption}

{
	\refstepcounter{wassumption}%
	\renewcommand{\thewassumption}{#1}%
	\begin{wassumption}\label{#2}%
	}
	{
	\end{wassumption}%
}

\usepackage{mathpazo}

\begin{document}
\renewcommand{\arraystretch}{1.5}

\title[Normal Hyperbolicity in Secondary Hopf Bifurcations]
{Normal Hyperbolicity in Secondary Hopf Bifurcations}

\author[D.D. Novaes and P.C.C.R. Pereira]
{Douglas D. Novaes$^1$  and Pedro C.C.R. Pereira$^2$}

\address{Departamento de Matem\'{a}tica - Instituto de Matemática, Estatística e Computação Científica (IMECC) - Universidade Estadual de Campinas (UNICAMP), Rua S\'{e}rgio Buarque de Holanda, 651, Cidade
Universit\'{a}ria Zeferino Vaz, 13083--859, Campinas, SP, Brazil}
\email{ddnovaes@unicamp.br$^1$}
\email{pedro.pereira@ime.unicamp.br$^2$}



\begin{abstract}
We combine results available in the literature to prove that the torus emerging in a secondary Hopf bifurcation is normally hyperbolic. This result is then applied to establish sufficient conditions for the bifurcation of normally hyperbolic invariant tori in the extended phase space of systems with small time-periodic perturbations via an application of the averaging method.
\end{abstract}

\maketitle

\section{Introduction}\label{sec:intro}
Establishing whether compact invariant manifolds exist in differential systems is crucial for understanding certain qualitative aspects of their dynamics. This understanding is further deepened when the asymptotic behavior of nearby trajectories is known, especially when the invariant manifolds are normally hyperbolic. In fact, this property guarantees not only their persistence with respect to small perturbations but also the existence of associated stable and unstable manifolds (see, e.g., \cite{Fenichel1971,hirschpughshub,wigginsnormally}).

In terms of dynamical properties, invariant tori in higher dimensional differential systems play a role particularly similar to that of limit cycles in planar ones, enabling the formulation of analogous questions regarding their existence and properties. For instance, in \cite{NP25}, an extension of Hilbert’s 16th Problem to higher dimensions was proposed, focusing on the investigation of the maximum number of isolated and normally hyperbolic invariant tori in 3D polynomial vector fields.  

Recently, results within the scope of averaging theory, developed in \cite{CANDIDO20204555, NP24, PNC23}, provided sufficient conditions for the existence of isolated invariant tori. The methods from \cite{NP24} were applied in \cite{NP25} as first steps toward understanding the aforementioned extension of Hilbert’s 16th Problem.

A classical approach for obtaining invariant tori in differential systems, as employed in \cite{CANDIDO20204555}, is through the \textit{secondary Hopf bifurcation}. Briefly explained, consider a generic one-parameter family of two-dimensional maps $ f(p,\mu)$ admitting a curve $p^*(\mu)$ such that:
\begin{itemize}
	\item[\footnotesize $\bullet$] $p^*(\mu)$ is a fixed point of $x \mapsto f(p,\mu)$ for each $\mu$ near $\mu^* \in \R$,
	\item[\footnotesize $\bullet$] the eigenvalues $\frac{\p f}{\p p}(p^*(\mu^*),\mu^*)$ lie on the unit circle and form a conjugated pair of distinct complex numbers.
\end{itemize}
Generically, such a family undergoes a Neimark-Sacker bifurcation, characterized by change of stability of the fixed point in conjunction with the emergence of an invariant closed curve from it (see, for instance, \cite{kuznetsov2023elements}).

When this bifurcation happens in a Poincaré map $p\mapsto \Pi(p,\mu),$ $p\in \Sigma$, associated to a periodic orbit of a family of differential equations $\dot x = F(x,\mu)$, $x\in\R^n$, the nascent invariant (under $\Pi$) closed curve in the Poincaré section $\Sigma$ corresponds to a torus in $\R^n$, itself invariant under the flow of $F$. This is usually called a secondary Hopf bifurcation.

In this paper, we essentially combine results available in the literature to ensure that the torus emerging in a secondary Hopf bifurcation is always normally hyperbolic. In short, normally hyperbolic invariant manifolds are a generalization of hyperbolic equilibria and limit cycles. It differs from the general notions of (asymptotic) stability and instability by not only characterizing the dynamics `transverse' to those manifolds, but also comparing it to the `tangential' dynamics, requiring the transversal behavior to dominate. An attracting normally hyperbolic invariant manifold, for instance, `flattens' neighborhoods of its points under forward action of the flow. Precise definitions can be found in \cite{hirschpughshub, wigginsnormally}.

We also use that finding to strengthen a recent result from \cite{CANDIDO20204555}, establishing sufficient conditions for the bifurcation of normally hyperbolic invariant tori in the extended phase space of systems with small time-periodic perturbations via an application of the averaging method. In particular, this allows us to state that the invariant tori detected in \cite{CANDIDO20204555,CNV20,CNS24,CV22,MC24} are all normally hyperbolic.

For illustration, we briefly discuss the case treated in \cite{CNV20}, which addresses the particularly elegant example of the three-dimensional Rössler system:
 $$(\dot x, \dot y, \dot z) = (-y-z,x+ay,bx-cz-xz),$$
where $a$, $b$, and $c$ are real parameters. By setting $a$, $b$, and $c$ as functions of a real parameter $\e$ for which $a(0)=c(0)=\overline{a} \in (-\sqrt{2},\sqrt{2})\setminus\{0\}$ and $b(0)=1$, the authors of that paper verify that, under suitable open conditions, an asymptotically stable invariant torus bifurcates from the trivial solution $x(t)=y(t)=z(t)=0$ as $\e$ crosses zero. This torus is obtained via a secondary Hopf bifurcation in the averaging framework developed in \cite{CANDIDO20204555}, so that our main result guarantees, \textit{a fortiori}, that it must also be normally hyperbolic. 

\section{Statement of the main theorem}
Consider the one-parameter family 
\begin{align}\label{eq:mainzero}
	\dot x = F(x,\mu), \quad (x,\mu) \in \Omega \times (-\mu_0,\mu_0)
\end{align}
of autonomous differential equations, where $\mu_0>0$, $\Omega$ is an open bounded subset of $\R^n$, $n\geq3$, and $F$ is a function of class $C^1$. 

Assume that, for $\mu=0$, \cref{eq:mainzero} has a $T_0$-periodic orbit, $T_0>0$, given by $t\mapsto\gamma_0(t)$. Consider a section $\Sigma$ transversal to $\gamma_0$ at a point $x_0$. The Implicit Function Theorem applied to the flow of \cref{eq:mainzero} ensures that a family of Poincaré maps $p\mapsto\Pi(p,\mu)$ is well defined on $\Sigma \cap U$, where $U$ is a neighbourhood of $x_0$.

Then, the premise that \cref{eq:mainzero} undergoes a secondary Hopf bifurcation can be stated as $\Pi$ undergoing a Neimark-Sacker bifurcation, i.e., we assume that
\begin{hypothesis}[S]
	There is a curve $p^*(\mu)$ such that, for all $\mu \in (-\mu_0,\mu_0)$, $\Pi(p^*(\mu),\mu) = p^*(\mu)$ and the eigenvalues $\{\lambda_1(\mu),\ldots,\lambda_{n-1}(\mu)\}$ of the derivative $\frac{\p \Pi}{\p p}(p^*(\mu),\mu)$ satisfy:
	\begin{enumerate}
		\item \label{item1hypothesisS}$\overline{\lambda_2}(0) = \lambda_1(0)$ and $|\lambda_1(0)| =1 $,
		\item \label{item2hypothesisS} $|\lambda_k(0)| \neq 1$ for $k \in \{3,\ldots,n-1\}$,
		\item $\left(\lambda_1(0)\right)^q \neq 1$ for $q \in \{1,2,3,4\}$,
		\item $d=\frac{d}{d\mu} \operatorname{Re}(\lambda_1(\mu)) \Big|_{\mu=0} \neq 0.$
	\end{enumerate}
\end{hypothesis}
From \cref{item1hypothesisS,item2hypothesisS}, we can assume that $\frac{\p \Pi}{\p p}(p^*(0),0)$ has, respectively $n_s$ and $n_u$ eigenvalues inside and outside the unit circle, with $n_s+n_u = n-3$.

The final element required for fully characterizing the families that undergo a secondary Hopf bifurcation is the so-called first Lyapunov coefficient $\ell_1$, which can be seen as quantifying non-degeneracy of non-linear terms of the system at $\mu=0$. For its calculation, the reader is referred to \cite[Equation (5.52) in Section 5.4.2]{kuznetsov2023elements}.


Before stating our main result, we restate, for completeness, the theorem describing the secondary Hopf bifurcation. In order to be as clear as possible, we provide the two cases depending on the sign of $\ell_1$ - usually known as subcritical and supercritical - separately. Its proof can be found in \cite[Sections 4.7 and 5.4]{kuznetsov2023elements}, and amounts to the fact that the family of Poincaré maps undergoes a Neimark-Sacker bifurcation.
\begin{theorem}\label{thm:maintheoremNS}
	If Hypothesis (S) is satisfied and $\ell_1>0$ (subcritical), there is a neighbourhood $V$ of $\gamma_0$ and $\mu_1 \in (0,\mu_0)$ such that, for each $\mu \in (-\mu_1,\mu_1)$, the following statements hold:
	\begin{enumerate}
		\item If $d\cdot \mu >0$, no invariant torus for \cref{eq:mainzero} exists in $V$ and $\gamma_\mu$ is a hyperbolic limit cycle with exactly $n_s$ characteristic multipliers inside the unit circle;
		\item If $\mu=0$, no invariant torus exists in $V$ and $\gamma_0$ is a non-hyperbolic periodic orbit with, respectively, $n_s$, $n_u$, and 2 characteristic multipliers inside, outside, and on the unit circle. The  directions corresponding to the multipliers on the unit circle are weakly repelling.
		\item If $d \cdot \mu<0$, there is a unique invariant torus for \cref{eq:mainzero} in $V$, which surrounds the hyperbolic limit cycle $\gamma_\mu$. Exactly $n_s+2$ characteristic multipliers of $\gamma_\mu$ are inside the unit circle.
	\end{enumerate}
	Similarly, if $\ell_1<0$ (supercritical), the following statements hold:
	\begin{enumerate}
		\item If $d\cdot \mu>0$, there is a unique invariant torus for \cref{eq:mainzero} in $V$, which surrounds the hyperbolic limit cycle $\gamma_\mu$. Exactly $n_s$ characteristic multipliers of $\gamma_\mu$ are inside the unit circle.
		\item If $\mu=0$, no invariant torus exists in $V$ and $\gamma_0$ is a non-hyperbolic periodic orbit with, respectively, $n_s$, $n_u$, and 2 characteristic multipliers inside, outside, and on the unit circle. The  directions corresponding to the multipliers on the unit circle are weakly attracting.
		\item If $d\cdot \mu <0$, no invariant torus for \cref{eq:mainzero} exists in $V$ and $\gamma_\mu$ is a hyperbolic limit cycle with exactly $n_s+2$ characteristic multipliers inside the unit circle;
	\end{enumerate}
\end{theorem}
Since it arises from a Neimark-Sacker bifurcation, the dynamics on the bifurcating invariant torus given by this theorem is also well understood. For instance, for a fixed value of the parameter, the Poincaré map generically has only a finite number of periodic orbits on the invariant curve corresponding to the torus. Also, as the parameter approaches its bifurcation value, infinitely many periodic orbits on the torus must emerge and collapse. We refer the reader to \cite[Section 7.3]{kuznetsov2023elements} for a more comprehensive treatment.

The main result of this paper can then be stated very succinctly as follows.

\begin{theorem}\label{thm:maintheorem}
	The invariant torus provided by Theorem \ref{thm:maintheoremNS} is normally hyperbolic. It has, exactly and respectively, $n_s$ and $n_s+1$ stable directions in the subcritical and supercritical case.
\end{theorem}

\section{Proof of the main result}

The proof is achieved in two steps. First, we notice that the invariant `circle' appearing in a Neimark-Sacker bifurcation is normally hyperbolic. This result was proven by Chaperon in \cite[Section 6]{chaperon2011generalised}. That allows us to state the following.
\begin{lemma}
	Under the assumptions of Theorem \ref{thm:maintheorem}, there is a neighbourhood $U$ of $p^*(0)$ and $\mu_1 \in (0,\mu_0)$ such that, for each $\mu \in (-\mu_1,\mu_1)$ satisfying $\ell_1\cdot d \cdot \mu<0$, $\Pi (\cdot,\mu)$ has a unique normally hyperbolic invariant closed curve $S_{\mu}$ in $\Sigma \cap U$, surrounding the fixed point $p^*(\mu)$.
\end{lemma}

The second step is essentially making use of the fact that, provided that the Poincaré map has a normally hyperbolic invariant closed curve, the flow has a normally hyperbolic invariant manifold given by its saturation. This follows from an application of a result due to Fenichel concerning `thin surface sections' (see \cite[Theorem 2.3]{Delshams2013} or \cite[Theorem 7]{fenichel1974}).

In fact, $S_\mu$ satisfies the following three properties defining such a thin surface of section (see \cite{fenichel1974}) for the flow $\Phi_t(x)$ of \cref{eq:mainzero}:
\begin{enumerate}
	\item  $F(x,\mu) \notin T_{x}S_{\mu}$, where $T_{x}S_{\mu}$ denotes the tangent space of $S_{\mu}$ at $x$.
	\item Each forward or backward orbit starting at a point in $S_{\mu}$ intersects $S_{\mu}$ again.
	\item $\tau(p) := \inf \left\{t>0 : \Phi_t(x) \in S_{\mu} \right\}$ defines a continuous function on $S_{\mu}$.
\end{enumerate}

Thus, since the saturation of $S_\mu$ forms an invariant torus $T_\mu$ intersecting $\Sigma$ transversely, we are able to conclude that.
\begin{lemma}
	$T_{\mu}$ is a normally hyperbolic invariant manifold for the flow of \cref{eq:mainzero}.
\end{lemma}


\section{Normally hyperbolic torus via averaging theory} \label{sec:averaging}

Consider a two-parameter family 
\begin{equation} \label{eq:mainone}
	\begin{aligned}
		&\dot x = \sum_{i=1}^k\e^i F_i(t,x,\mu) + \e^{k+1} \tilde{F}(t,x,\mu,\e), \\ & (t,x,\mu,\e) \in \R \times D \times (-\mu_0,\mu_0) \times (-\e_0,\e_0) 
	\end{aligned}
\end{equation}
of non-autonomous differential equations, where $k \in \mathbb{N^*}$, $\mu_0,\e_0>0$, and $D$ is an open bounded subset of the plane. It is assumed that $F_i$, $i \in \{1,\ldots,k\}$, and $\tilde{F}$ are $C^1$ functions which are $T$-periodic in $t$ with $T>0$. Periodicity in $t$ allows us to regard \cref{eq:mainone} as a family of autonomous systems defined in the `extended phase space' $\R/(T\mathbb{Z}) \times D$ and given by 
\begin{equation} \label{eq:maintwo}
	\begin{aligned}
		& \dot t =1, \\
		&\dot x = \sum_{i=1}^k\e^i F_i(t,x,\mu) + \e^{k+1} \tilde{F}(t,x,\mu,\e).
	\end{aligned}
\end{equation} 
In our context, $\e$ is a perturbative parameter often assumed to have a fixed value, and we wish to analyze the behaviour under changes of $\mu$.

Cândido and Novaes (see \cite[Theorem B]{CANDIDO20204555}) have proven that, under the assumption that the first non-vanishing averaging function associated to \cref{eq:mainone} has an equilibrium with a pair of non-zero complex conjugated eigenvalues crossing the imaginary line transversely as $\mu=0$, combined with a generic non-degeneracy condition concerning the Lyapunov coefficient, \cref{eq:maintwo} undergoes the birth of an invariant torus.

The proof of this theorem is achieved by showing that, for small value of $\e>0$, the stroboscopic (i.e., defined on the section $\Sigma_0:=\{\overline{0}\} \times D$) Poincaré map $\Pi_\e: x \mapsto \varphi_\e(T,x)$, where $\varphi_\e$ is the solution of \cref{eq:mainone} with $\varphi_\e(0,x)=x$, undergoes a Neimark-Sacker bifurcation as $\mu$ passes through a critical value $\mu^*(\e)$. In particular, this allows the authors in \cite{CANDIDO20204555} to state the torus obtained is either asymptotically stable or unstable.

As a consequence of the result proven in this paper, we are able to strengthen this result by showing that the torus appearing in the family is a normally hyperbolic invariant manifold for \cref{eq:maintwo}.

The hypotheses of Theorem B in \cite{CANDIDO20204555} are stated in terms of averaging functions. Hence, to properly state our strengthened result, we first briefly introduce them  along with necessary concepts from the theory of averaging.

For small $\e$, the stroboscopic Poincaré map $\Pi_\e$ of \cref{eq:maintwo} defined on the section $\Sigma_0 =\{\overline{0}\} \times D$ can be expanded in $\e$ as follows (see \cite[Lemma 1]{Llibre_2014}).
\begin{align}
	\Pi_\e(x,\mu) = x + \sum_{i=1}^k \e^i g_i (x,\mu) + \e^{k+1} R(x,\mu,\e),
\end{align}
where each $g_i$, labeled \emph{averaged function of order $i$}, is calculated using the formulas available in \cite[Appendix A]{CANDIDO20204555}, and $R$ is a remainder function (see also \cite{Novaes21b}). The formulas will be omitted here for brevity.

Assume that not all $g_1, \ldots, g_k$ vanish identically and let $l \in \{1,2,\ldots, k\}$ be the index of the first non-vanishing averaged function. The first hypothesis required for our theorem is that $\dot x = g_l(x,0)$ has a `Hopf point' at $x=x_0$, i.e., 

\begin{hypothesis}[H]
	$g_l(x_0,0)=0$ and $\frac{\p g_l}{\p x}(x_0,0)$ has a pair of conjugated purely imaginary eigenvalues $\pm i \omega_0$, $\omega_0>0$.
\end{hypothesis}

Assuming that (H) holds, the Implicit Function Theorem yields a curve $\xi(\mu)$ defined on a neighbourhood of $\mu=0$ such that $g_l(\xi(\mu),\mu) =0 $. Moreover, by continuity of the eigenvalues of a matrix with respect to its entries, $\frac{\p g_l}{\p x}(\xi(\mu),\mu)$ must have a pair of conjugated complex eigenvalues $\alpha(\mu) \pm i \beta(\mu)$ for small values of $\mu$. 

The second hypothesis regards transversality of the family $g_l(x,\mu)$ with respect to the `surface' of Hopf germs.

\begin{hypothesis}[T] \label{hypothesish} $\alpha'(0) \neq 0$. 
\end{hypothesis}

Assuming that (H) and (T) hold and that $\e$ is sufficiently small, Cândido and Novaes showed that there is a locally unique fixed point $x^*(\mu,\e)$ of $x \mapsto \Pi_\e(x,\mu)$, which satisfies $x^*(\mu,0) =\xi(\mu)$.  More importantly, they proved that there is $\mu_\e$, depending continuously on the parameter $\e$, such that $\Pi_\e$ has a `transversal Neimark-Sacker' point at $(x^*(\mu_\e,\e),\mu_\e)$, i.e., $\Pi_\e(x^*(\mu_\e,\e),\mu_\e) = x^*(\mu_\e,\e)$ and the pair $\{\lambda_\e(\mu), \overline{\lambda_\e}(\mu)\}$ of complex conjugated eigenvalues of $\frac{\p \Pi_\e}{\p x}(x^*(\mu,\e),\mu)$ satisfies
\begin{itemize}
	\item[\footnotesize $\bullet$] $|\lambda_\e(\mu_\e)|=1$,
	\item[\footnotesize $\bullet$] $(\lambda_\e(\mu_\e))^q \neq 1$ for $q \in \{1,2,3,4\}$,
	\item[\footnotesize $\bullet$] $\frac{d}{d\mu} |\lambda_\e|\Big|_{\mu=\mu_\e} \neq 0$
\end{itemize}

The third and final hypothesis is concerned with non-degeneracy of the first Lyapunov coefficient $\ell^\e_1$ associated to the map $\Pi_\e$ at the bifurcation point $(x^*(\mu_\e,\e),\mu_\e)$. This quantity can also be expanded near $\e=0$ as
\begin{align}
	\ell_1^\e = \e^l \ell_{1,l} + \e^{l+1} \ell_{1,l+1} + \cdots + \e^k \ell_{1,k} + \mathcal{O}(\e^{k+1}).
\end{align}
Expressions for each $\ell_{1,i}$, which depends only on the $g_j$ for $j \in \{1,\ldots,i\}$, are provided in \cite[Section 4.3]{CANDIDO20204555}. The non-degeneracy condition required of $\ell_1^\e$ is given by
\begin{hypothesis}[ND] \label{hypothesisnd}
	There is $j \in \{l,\ldots, k\}$ such that $\ell_{1,j} \neq 0$. 
\end{hypothesis}
Let $j^*$ be the first index $j$ for which $\ell_{1,j^*} \neq 0$. In \cite{CANDIDO20204555}, as a consequence of \Cref{thm:maintheoremNS}, Cândido and Novaes have proved the following :

\begin{theorem}\label{thm:theoremaveraging}
	Assuming the validity of Hypotheses (H), (T), and (ND), there is, for each small $\e>0$, an interval $I_\e\subset (-\mu_0,\mu_0)$ containing $\mu(\e)$, and an open set $U_\e \subset \R / (T\mathbb{Z}) \times D$ such that
	\begin{enumerate}
		\item If $\mu\in I_{\e}$ and $\alpha'(0)\cdot \ell_{1,j^*} (\mu - \mu(\e)) \geq 0$, then \cref{eq:maintwo} has a $T$-periodic orbit $\gamma_{\mu,\e} \subset U_\e$ corresponding to the fixed point $x^*(\mu,\e)$ of $x \mapsto \Pi_\e(x,\mu)$ and admits no invariant tori in $U_\e$;
		\item If $\mu\in I_{\e}$ and $\alpha'(0)\cdot\ell_{1,j^*} (\mu - \mu(\e)) < 0$,  then \cref{eq:maintwo} has a $T$-periodic orbit $\gamma_{\mu,\e} \subset U_\e$ corresponding to the fixed point $x^*(\mu,\e)$ of $x \mapsto \Pi_\e(x,\mu)$ and a unique invariant torus $T_{\mu,\e}$ in $U_\e$, which surrounds $\gamma_{\mu,\e}$. 
	\end{enumerate}
\end{theorem}

With the introduction of \Cref{thm:maintheorem}, more can be said. To wit, the following stronger result holds true.
\begin{theorem}
	The torus appearing in \Cref{thm:theoremaveraging} is normally hyperbolic and repelling, if $\ell_{1,j^*}>0$, or attracting, if $\ell_{1,j^*}<0$.
\end{theorem}

In fact, the emergence of the torus $T_{\mu,\e}$ from the $T$-periodic orbit $\gamma_{\mu,\e}$ described in \Cref{thm:theoremaveraging} is a direct consequence of \Cref{thm:maintheorem} once we take into account the following result (see \cite[Theorem B]{CANDIDO20204555} for a proof).
\begin{theorem}
	Assume that Hypotheses (H), (T), and (ND) hold. For small values of $\e$, the Poincaré map $\Pi_\e$ undergoes a Neimark-Sacker bifurcation at $x^*(\mu_\e,\e)$ for $\mu=\mu_\e$.
\end{theorem}

The significance of Theorem \ref{thm:theoremaveraging} in the study of normally hyperbolic limit cycles, compared to those presented in \cite{PNC23},  lies in the fact that it permits the simultaneous use of multiple averaging functions.  In contrast, the results from \cite{PNC23} impose conditions on the guiding system only, relying solely on the first non-vanishing averaged function. Consequently, in some cases, Theorem \ref{thm:theoremaveraging} allows for the detection of normally hyperbolic invariant tori that cannot be detected using the results in \cite{PNC23}, making it a valuable addition to the `averaging toolbox' for studying invariant tori.

\section*{Acknowledgements}
PCCRP is supported by S\~{a}o Paulo Research Foundation (FAPESP) grants nº 2023/11002-6 and 2020/14232-4. DDN is partially supported by S\~{a}o Paulo Research Foundation (FAPESP) grants nº 2018/13481-0 and by Conselho Nacional de Desenvolvimento Cient\'{i}fico e Tecnol\'{o}gico (CNPq) grant nº 309110/2021-1.

\bibliography{references}
\bibliographystyle{plain}

\end{document}